# Intro To Logarithms After The Meltdown

Dr. David A. de Wolf
103 Sherwood Ct., Blacksburg, VA 24060
Tel.: 540-951-3073, Email: dadewolf@vt.edu

*Abstract:* This work is an attempt to 'reconstruct' logarithms in the hypothetical case that mankind has suffered a catastrophe through which all repositories of (mathematical as well as other) knowledge are lost, with the exception of simple arithmetic operations. It turns out that this is possible without re-inventing calculus!

Suppose, through some self-inflicted catastrophe, mankind has been greatly reduced, and all its repositories of previous knowledge (libraries, tapes, CD's, other electronic storage devices,...) destroyed. Man slowly learns to regain His mathematical knowledge, but of course all recourse to previous results is lost. I imagine that one of the first ideas to arise after relearning the basic operations is that giving rise to what we call 'the logarithm'. Here is an attempt to define logarithms as I envision that our unhappy descendents might rediscover it. Remember, they would have no prior memory or knowledge of the exponential function, or of all the beautiful algorithms with which to calculate logarithms. Instead, they would probably rediscover logarithms by noting something special about multiples of 10. But before going in to this, I think it might be safe to assume that memory of our ten-based counting system (with its prominent zero) and of the simplest operations (addition, subtraction, multiplication, division) are retained or rediscovered. I also will assume that the operation of taking square roots will be rediscovered rapidly (I show below a simple algorithm requiring no more than iteration of the simplest arithmetic operations). As for notation: sometimes I show a product of $a \times b$ as $ab$ and sometimes as $(a)(b)$ where $a$ and $b$ stand for numerical quantities without specifying exactly which ones, because we use these 'names' in more general relationships where the relationship holds no matter what $a$, $b$ really are. Sometimes it is useful to use a letter with a subscript to indicate a 'name' of a quantity, e.g. as in the following

Consider then the series of multiples, $a_0 = 1$, $a_1 = 10$, $a_2 = 100$, $a_3 = 1,000$, $a_4 = 10,000$. $a_5 = 100,000$,....................

The survivors would note that there are two ways to look at $a_n$; $(i)$: $a_n$ consists of a 1 with $n$ zeroes following it, $(ii)$: $a_n$ is equal to multiplying 10 by itself $n$ times. For example, $a_3 = 1,000 = 10 \times 10 \times 10$. They might use this to write any of these as the equivalent of $10^n$ (*i.e.* for a '1' with $n$ zeroes following it). Their notation might not be identical to $10^n$ (*e.g.* they might write it as $[1, n]$) but it would almost certainly amount to the same thing so I will adopt our notation for simplicity. Thus $1 = 10^0$ because this 1 has no zeroes following it, $1000 = 10^3$ because there are three zeroes following 1 here. The raised number (we call this 'an exponent') $n$ is defined to be '**log($a_n$)**'. Again, the survivors might have some other notation to indicate the exponent, but we'll use the familiar one as it would amount to the same thing. Thus $\log(a_3) = \log(1,000) = \log(10^3) = 3$, $\log(a_4) = \log(10,000) = \log(10^4) = 4$, $\log(a_5) = \log(100,000) = \log(10^5) = 5$, .... and in general

$$\boxed{\log(a_n) = \log(10^n) = n} \tag{1}$$

The significant feature of $\log a_n$ is that for $a_n = 10^n$ the $\log(a_n)$ exponents are equally spaced although the $a_n$ are not (see sketch; negative exponents are explained below). In fact, this would allow us to define $10^0 = 1$ (this is a '1' with 0 zeroes behind it) with $\log(10^0) = 0$. This renders it easier to express large numbers. For example in our 'decibel' scale where 110 dB describes a quantity that is $10^{11}$ times as large as a reference quantity of size 1 (for some applications there is a factor 2 implied so that the magnification factor is 'only' $10^{5.5}$).

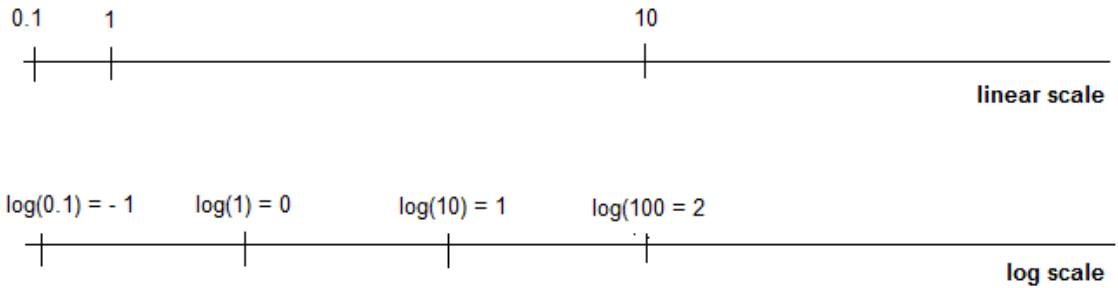

It will be apparent that $a_n \times a_m = a_{n+m}$. I often write a product such as $a_n \times a_m$ as $a_n a_m$ (just leaving off the $\times$ multiplication sign). I'll use both notations below.

***Example***: for $a_3 = 10^3$ and $a_2 = 10^2$:  $10^3 \times 10^2 = 1,000 \times 100 = 100,000 = 10^5$ -- just add the zeroes. In equations like this, each number on either side of the $\times$ sign is called a *factor*.

It follows easily that $\log(a_n a_m) = \log(a_{m+n}) = \log(10^{m+n}) = m + n = \log(a_n) + \log(a_m)$.

Note also that $(a_n)^m = [a_n \times a_n \times ... \times a_n]_{m\,factors}$. As shown how above, apply $a_n \times a_n = a_{2n}$, then $a_{2n} \times a_n = a_{3n}$ and so on until you ultimately get $a_{mn} \equiv a_{nm} = 10^{nm}$. Therefore $\log(a_n)^m = \log(a_{nm}) = \log(10^{nm}) = nm$. The same result comes out for $m\log(a_n) = m\log(10^n) = mn$ upon applying Eq. (1). In summary:

$$\boxed{\log(a_n a_m) = \log(a_n) + \log(a_m) = m + n} \tag{2a}$$

$$\boxed{\log[(a_n)^m] = m\log(a_n) = nm} \tag{2b}$$

Note that also $\log[(a_m)^n] = m\log(a_n) = mn$ (same result) as you can see by interchanging $m$ and $n$.

***Example*** :
$\log[(a_2)^3] = \log(100 \times 100 \times 100) = \log(1,000,000) = \log(10^6) = 6.$   This is the same as $\log(a_3^2) = \log(1000 \times 1000) = \log(1,000,000) = \log(10^6) = 6.$

The most special thing about Eqs. (2) is that *the logarithm of a product is the sum of the individual logarithms.*

The next obvious question is whether our survivors can extend this definition to noninteger $n$, *e.g.* $n = 1/2$. Let's see what is needed. From Eq. (2) for $m = n$ it should be apparent that

$$\log[(a_n)^2] = 2\log(a_n) = 2n$$

Apply this "without thinking" to $n = 1/2$ so that $2n = 1$. This means we have defined $a_{1/2}$ as $a_{1/2} \equiv 10^{1/2}$ but now we have to find out what this means.

We already know that $1 = \log(10)$ so that we seem to have found that

$$\log[(a_{1/2})^2] = 2 \times \tfrac{1}{2} = 1 = \log(10) \rightarrow (a_{1/2})^2 = 10$$

from which we conclude that $a_{1/2} = \sqrt{10}$, so that an obvious extension would be to have

$$\boxed{a_{1/2} \equiv 10^{1/2} = \sqrt{10} \approx 3.16} \tag{3a}$$

It also follows that $\log(a_{1/2)} \equiv \log(10^{1/2}) = \tfrac{1}{2}\log(10) = 1/2$.

$$\boxed{\log(10^{1/2}) = \log(\sqrt{10}) = 1/2} \tag{3b}$$

So the art of taking square roots is also needed. It is not unreasonable to assume that our survivors have already relearned to take square roots. They might have already found the fairly straight-forward and simple repetitive algorithm:

*i)* Guess what the square root of number $x$ is. Call it $(x_{\sqrt{}})_1$.
*ii)* Calculate $(y_{\sqrt{}})_1 = x/(x_{\sqrt{}})_1$.
*iii)* Calculate $(x_{\sqrt{}})_2 = \tfrac{1}{2}[(x_{\sqrt{}})_1 + (y_{\sqrt{}})_1]$.
*iv)* Calculate $(y_{\sqrt{}})_2 = x/(x_{\sqrt{}})_2$.
*v)* Calculate $(x_{\sqrt{}})_3 = \tfrac{1}{2}[(x_{\sqrt{}})_2 + (y_{\sqrt{}})_2]$.
*vi)* Keep going...you add on two or more significant figures each further step! For example, the square root of 1747 is found accurately to 10 significant figures after only 3 iterations beyond the initial guess of 40.

***Example***: Do $\sqrt{1747}$ as follows:
Let $(x_{\sqrt{}})_1 = 40$
Then $(y_{\sqrt{}})_1 = 1747/40 = 43.675$,
$(x_{\sqrt{}})_2 = \tfrac{1}{2}(40 + 43.675) = 41.8735$,
$(y_{\sqrt{}})_2 = 1747/41.8735 \approx 41.75679713$
$(x_{\sqrt{}})_3 = \tfrac{1}{2}(41.8735 + 41.75679713) \approx 41.79714857$

$$(y_\sqrt{})_3 = 1747/41.79714857 \approx 41.79710961$$
$$(x_\sqrt{})_4 = \tfrac{1}{2}(41.79714857 + 41.79710961) \approx 41.79712909$$
$$(y_\sqrt{})_4 = 1747/41.79712909 \approx 41.79712909$$
$$(x_\sqrt{})_5 = \tfrac{1}{2}(41.79712909 + 41.79712909) \approx 41.79712909$$

and indeed, $(x_\sqrt{})_5 \approx (x_\sqrt{})_4$ to 10 significant figures!

In the next step I again apply

$$\log[(a_n)^2] = 2\log(a_n) = 2n$$

but this time we choose $n = 1/4$ so that we get $\log[(a_{1/4})^2] = 1/2 = \log(\sqrt{10})$ from Eq. (3b). By this means we conclude that $(a_{1/4})^2 = \sqrt{10}$ and therefore that $a_{1/4} = \sqrt{\sqrt{10}} = \sqrt[4]{10} \approx 3.16227766$. By its very definition $a_{1/4} \equiv 10^{1/4}$.

Continuing this procedure, we find that

$$a_{1/8} = 10^{1/8} = \sqrt[8]{10} \approx 1.333512432,$$
$$a_{1/16} = 10^{1/16} \equiv \sqrt[16]{10} \approx 1.154781985, ....$$

and ultimately we find for $a_{1/N}$ where $N = 2^n$ and $n$ is a very large integer so that $2^N$ is ever so much larger (***Example***: $n = 12 \rightarrow N = 4,096$) :

$$\boxed{a_{1/N} \equiv 10^{1/N} = \sqrt[N]{10} \text{ where } N = 2^n} \tag{4}$$

Thus this $a_{1/N}$ approaches 1 as $n$ is increased (*This can be proved more rigorously but seems to be clearly the case upon calculation*). The following should be fairly apparent.

$$(a_{1/N})^n = (10^{1/N} \times 10^{1/N} .... \times 10^{1/N})_{n\,factors} = 10^{(\frac{1}{N}+...+\frac{1}{N})_{n\,terms}} = 10^{n/N}$$

and if $n$ is sufficient large, so that $1/N = 1/2^n \ll 1$ is essentially 'infinitesimal', then $n/N \equiv x$, even for $n > N$, will approximate as closely as we wish any real number $> 1$. (*By 'closely' I mean accurate to a number of digits considerably less than those in $1/N$ but that is not restrictive because N can be made arbitrarily large* ). That gives meaning to

$$\boxed{y = 10^x \text{ for } x > 1} \text{ and to } \boxed{\log(y) = \log(10^x) = \log(10^{n/N}) = n/N = x} \tag{5}$$

The equations (2) now also hold for $x_1$, $x_2$ and $y_1$, $y_2$ instead of just for the integers $m, n$ and $10^m$, $10^n$. Specifically if $y_1 = 10^{x_1}$ and $y_2 = 10^{x_2}$ then $\log(y_1 y_2) = \log[(10^{x_1})(10^{x_2})] = \log(10^{x_1+x_2}) = x_1 + x_2$:

$$\boxed{\log(y_1 y_2) = x_1 + x_2 = \log(y_1) + \log(y_2)} \tag{6}$$

This is the extension of Eq. (2) but now for all $y \geq 1$. We shall need this 2nd property very much. It merely expresses that *the logarithm of a product of rational numbers $> 1$ is the sum of the logarithms of each of those rational numbers*.

Here are two more relationships of importance, derived from $\log(y) = x$ and $10^x = y$:

$$10^{\log(y)} = 10^x = y \tag{7a}$$
$$\log(10^x) = \log(y) = x \tag{7b}$$

Thus log and taking a power of 10 are inverse operations; one 'undoes' the other.

Now let us consider what to do about the numbers $y = 10^x$ with $0 < y < 1$. We are going to find that $x$ is *negative* and specifically that $-\infty < x < 0$.

Here's how: suppose we have a product $y_1 y_2 = 1$ with $y_1 = 10^{x_1}$, $y_2 = 10^{x_2}$ where $x_1 > 1$ is a rational number. We then know that $\log(y_1 y_2) = \log(1) = 0$ or $y_1 y_2 = 1$. Therefore it is also true that $\log(y_1 y_2) = \log(10^{x_1} 10^{x_2}) = \log(10^{x_1+x_2}) = 0 = \log(1) = \log(10^0)$. From this we are forced to conclude, for consistency in notation, that $x_2 = -x_1$ and thus that $y_2 = 10^{x_2} = 10^{-x_1}$. We have then found that the logarithms of numbers between 0 and 1 span the range from $-\infty$ to 0 whereas those of numbers between 1 and $\infty$ span the range from 0 to $\infty$. So now the entire range of real numbers is covered!

The logarithms of negative numbers require more sophisticated mathematics (complex numbers), the reinvention of which I shall not deal with here.

Here is a practical way to calculate the logarithm in the interval [1,10] starting from no other knowledge than multiplication and the art of taking square roots (this method would be preferable to spanning together enormously long products of $a_{1/N}$).

> First calculate the logarithms of 10 and $\sqrt{10}$ with $\log(10) = 1$ and $\log(\sqrt{10} = 1/2$ then of $10^{1/4}$ and $10^{3/4} = (\sqrt{10})(\sqrt[4]{10})$ with logs equal, respectively, to 1/4 and 3/4, then of $(\sqrt[8]{10}) = 10^{1/8}, 10^{3/8} = (\sqrt[4]{10})(\sqrt[8]{10}), 10^{5/8} = (\sqrt{10})(\sqrt[8]{10})$, and $10^{7/8} = (\sqrt{10})(10^{3/8})$ with logs equal to 1/8, 3/8, 5/8, and 7/8.

So then we have the logs 1/8, 1/4, 3/8, 1/2, 5/8, 3/4, 7/8, and 1 ( spaced at $1/8$). As you see, each steps halves the spacing between successive logarithms calculated up to that point. We thus can fill in successively smaller gaps between 1 and 10. Ultimately, we arrive at a factor arbitrarily close to unity by calculating $\sqrt[N]{10}$ for $N = 2^n$. and thus using $1/N$ as the spacing for the logarithm of numbers in the [1, 10] interval. The intervals [10,100], [100,1000],..., and also [0,1] can be handled in a similar fashion using (2).

Note that we are not limited to using 10 as 'base' number. Every mathematical symbol and equation above with a '10' in it is equally valid if we replace '10' by an arbitrary real number $p$. Let me explain:

The number 134 can be written as $4 + 3 \times 10 + 1 \times 10^2$.
The number 54.79 can be written as $9 \times 10^{-2} + 7 \times 10^{-1} + 4 \times 10^0 + 5 \times 10^1$

Likewise, any number (with or without decimal parts) can be written as

$$... + (a_{-3})(10^{-3}) + (a_{-2})(10^{-2}) + (a_{-1})(10^{-1}) + a_o + (a_1)(10) + (a_2)(10^2) + ...$$

Now suppose we chose 3 as a base. Then we would label the integers (starting from one) as follows (top row only):

$$\begin{array}{cccccccccccccccc} 1 & 2 & 10 & 11 & 12 & 20 & 21 & 22 & 100 & 101 & 102 & 110 & 111 & 112 & 120 & 121 & ... \\ 1 & 2 & 3 & 4 & 5 & 6 & 7 & 8 & 9 & 10 & 11 & 12 & 13 & 14 & 15 & 16 & ... \end{array} \quad (8)$$

The bottom row gives the 10-based equivalents. The number we call 3 in the 10-based system would be '10' in the 3-based system. Likewise '100' in the top row would be $9 = 3 \times 3 = 3^2$ (in an obvious extension of what we did with 10-based products of 10 with itself) and '1000' would be $27 = 3 \times 3 \times 3 = 3^3$. Also, '11' in the top row would be 4 in the bottom row, '12' would be 5, '20' would be 6, and '21' would be 7. Perhaps it is now clear how this works. Referring to (8), I can write '22' and '120' (these are <u>not</u> the 22 and 120 as we know them in the 10-based system) in the 3-based system as

'22' = $2 \times 3^0 + 2 \times 3$ ( = 8 in the 10-based system)
'120' = $0 \times 3^0 + 2 \times 3 + 1 \times 3^2$ ( = 15 in the 10-based system)

In general any integer can be written in the 3-based system as

$$... + (a_{-3})(3^{-3}) + (a_{-2})(3^{-2} + (a_{-1})(3^{-1}) + a_o + (a_1)(3) + (a_2)(3^2) + ... \quad (9)$$

and likewise, in the $p$-base system, as

$$... + (a_{-3})(p^{-3}) + (a_{-2})(p^{-2} + (a_{-1})(p^{-1}) + a_o + (a_1)(p) + (a_2)(p^2) + ... \quad (10)$$

Thus we see there is nothing special about base 10; in all of the above we could have substituted 3 or $p$ for 10. The big new step is that, instead of writing $y = 10^x$ and $x = \log(y)$ we now write the same $y$ as $y = p^z$ and define $z \equiv \log_p(y)$ for the logarithm in base $p$.

From now on, we will write $\log_{10}(y)$ for clarity instead of just $\log(y)$. Then it follows that

$p = (p^z)^{1/z} = (10^x)^{1/z} = 10^{x/z}$, hence it follows that
$\log_{10}(p) = x/z \rightarrow z = x/\log_{10}(p) = \log_{10}(y)/\log_{10}(p)$.

The steps in these two lines are the same as if $p$ were 10. Thus we have found $z = \log_{10}(y)/\log_{10}(p)$. Remember that $z \equiv \log_p(y)$, then these two identities give

$$\log_p(y) = \frac{\log_{10}(y)}{\log_{10}(p)} \tag{11}$$

thus explaining how we convert a logarithm in base 10 to base $p$. Note that things come out right in (11) if $p = 10$ because $\log_{10}(10) = 1$. There is nothing special about '10' here [see sentence under Eq. (10)], and if you substitute some other base, e.g. $q$ for 10 in the previous ten or so lines you would get

$$\log_p(y) = \frac{\log_q(y)}{\log_q(p)} \tag{12}$$

to convert from base $q$ to base $p$. All the other relationships hold, e.g.

$$\log_p(y_1 y_2) = \log_p(y_1) + \log_p(y_2) \tag{13}$$

But what about $\log_p(y^x)$ for arbitrary $x, y, p$? In Eq. (5) we saw that $x = \log_{10}(y)$ if $y = 10^x$ for any rational numbers $x, y$. We have now seen that we can substitute any base $p$ for 10 and therefore we also have

$$x = \log_p(y) \text{ if } y = p^x \tag{14}$$

Now I want to approach the subject of a very special base $p$ that we shall call $e$. To do so we need to look at *graphic plots*. A graph in its simplest form of $y = \log_p(x)$ versus $x$ will show values of $x$ on a horizontal axis and of $y$ on a vertical one. Based on what we have seen so far, it might resemble the following sketch:

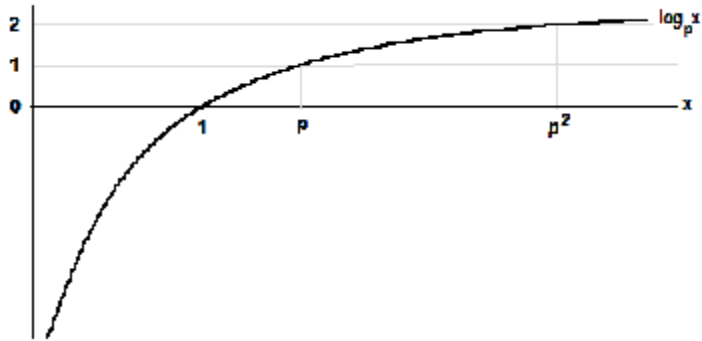

At $x = 1$ we would find $y = \log_p(1) = 0$ because $p^0 = 1$. At $x = p$ we would find $y = \log_p(p) = 1$ because $p^1 = p$. Then $\log_p(p^2) = 2$. And so on. It seems clear that the 'slope' of the curve at $x = 1$ will be steeper for smaller $p$ (because $p$ will be closer to 1 on the horizontal axis whereas the spacing on the vertical axis is unchanged). The slope therefore tells us something about $p$.

If this curve were a straight line (which it isn't) then there would be only one slope which we could characterize by the ratio $(y_2 - y_1)/(x_2 - x_1)$ for any two points $(x_1, y_1)$ and $(x_2, y_2)$ on the line. If $x$ is the time and $y$ is the distance traveled by a car then $(y_2 - y_1)/(x_2 - x_1)$ would be the speed of the car, so slope does give useful information.

As the curve is not a straight line, we find the slope at a desired point by drawing a tangent line there (see next sketch) and then the slope is characterized by that same ratio $(y_2 - y_1)/(x_2 - x_1)$ except that both points are 'infinitesimally' close to the point where the tangent line touches the curve. We could also choose $x_1 = x$ and $x_2 = x + \epsilon$ where $\epsilon$ is an infinitesimally small number (I'll show how to choose it shortly).

For $p = 10$ the slope at $x$ can be found from $y_2 - y_1 = \log_{10}(x + \epsilon) - \log_{10}(x)$ $= \log_{10}[(x)(1 + \frac{\epsilon}{x})] - \log_{10}(x) = [\log_{10}(x) + \log_{10}(1 + \frac{\epsilon}{x})] - \log_{10}(x)$ so that

$$y_2 - y_1 = \log_{10}(x + \epsilon) - \log_{10}(x) = \log_{10}(1 + \epsilon/x) \tag{15}$$

We shall use this relationship for infinitesimally small $\epsilon$. As you can see from the sketch below, the slope is determined by the angle that a tangent to the curve at $x$ makes with the horizontal axis.

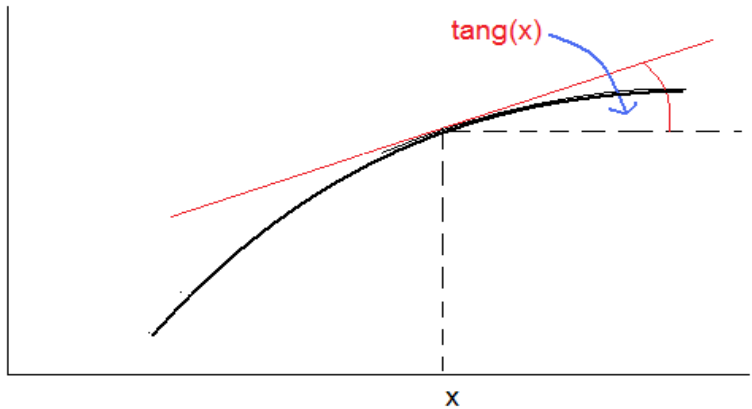

Rediscovery of this would not be a big stretch for our survivors. The tangent of that angle would be (I used Eq. (15) to simplify this in the 2nd step):

$$\text{tang}_{10}(x) = \frac{y_2 - y_1}{x_2 - x_1} = \frac{\log_{10}(x+\epsilon) - \log_{10}(x)}{(x+\epsilon) - x} = \frac{\log_{10}(1+\epsilon/x)}{\epsilon} \tag{16}$$

in the limit as $\epsilon \to 0$. We recognize this as the *derivative* but I'll assume calculus has not yet been rediscovered.

$$\boxed{\text{tang}_{10}(x) = \lim_{\epsilon \to 0}\left[\frac{\log_{10}(1+\epsilon/x)}{\epsilon}\right]} \tag{17}$$

This notation means "*First do what is behind the 'lim' and then let $\epsilon \to 0$*". The quantity $(1 + \epsilon/x)$ in (17) is a number very close to 1 for very small $\epsilon$. As I showed above, so is $\sqrt[N]{10}$ for large $N$. Therefore, set $1 + \epsilon/x = \sqrt[N]{10}$ with $N = 2^n$ and with large $n$.

E.g. choose $n = 20$ so that $N = 2^{20} = 1,048,576 \rightarrow 1 + \epsilon/x = \sqrt[N]{10}$ $\approx 1.000002196$ so that $\epsilon/x \approx 0.000002196$. It then also follows that $\log_{10}(1 + \epsilon/x) = \log_{10}(\sqrt[N]{10}) = \log_{10}(10^{1/N}) = 1/N = 1/2^{20}$ and therefore that the tangent is estimated (at least to this accuracy) by

$$\text{tang}_{10}(x) \approx \left[\frac{\log_{10}(1+\epsilon/x)}{\epsilon}\right] = \left[\frac{1}{x}\frac{\log_{10}(1+\epsilon/x)}{\epsilon/x}\right] \approx \frac{1}{x}\frac{(1/2^{20})}{0.000002196} \approx \frac{1}{x}0.434$$

Using our constructed logarithm table, we find that $0.434 \approx \log_{10}(2.718...)$ We give this new number, $2.718...$, the name '$e$'; greater accuracy for $e$ can be obtained by increasing $n$. We of course can use calculus to show that (17) converges, but our survivors will find this 'experimentally'. As a result:

$$\boxed{\text{tang}_{10}(x) = \lim_{\epsilon \to 0}\left[\frac{\log_{10}(1+\epsilon/x)}{\epsilon}\right] \equiv \frac{1}{x}\log_{10}(e) \approx 0.434\frac{1}{x}, \quad e \approx 2.718...} \tag{18}$$

Just so that the ghosts of the vanished generations are satisfied, it is relatively easy to see that this limit exists. It is basically $t = \frac{1}{N(10^{1/N}-1)}$. Note that $\ln(10^{1/N}) = \frac{1}{N}\ln(10) \rightarrow 10^{1/N} = e^{\ln(10)/N} = 1 + \frac{\ln 10}{N} + \frac{(\ln 10)^2}{2N^2} + .....$ so that $10^{1/N} - 1 = \frac{\ln 10}{N} + O(\frac{1}{N^2})$ and thus $t = \frac{1}{\ln(10)} + O(\frac{1}{N}) = \log_{10}(e)$ in the limit $N \rightarrow \infty$.

For base $p$ the slope at $x$ is related to that for base 10 by (8), i.e.
$$\text{tang}_p(x) = \lim_{\epsilon \to 0}\left[\frac{\log_p(1+\epsilon/x)}{\epsilon}\right] = \lim_{\epsilon \to 0}\left[\frac{\frac{1}{\epsilon}\log_{10}(1+\epsilon/x)}{\log_{10}(p)}\right] = \frac{1}{x}\frac{\log_{10}(e)}{\log_{10}(p)} \text{ so that}$$

$$\boxed{\text{tang}_p(x) = \frac{1}{x}\frac{\log_{10}(e)}{\log_{10}(p)} = \frac{1}{x}\log_p(e)} \tag{19}$$

Hence what is special about $p = e$ is that the tangent of the curve $\log_p(x)$ vs. $x$ is $1/x$ when the base is chosen to be $p = e$. The survivors realize that $p = e$ is a special base and they decide to give the logarithm-to-that-base a special name: $\boxed{\log_p \equiv \ln}$. What is special about $\ln(x)$ with $e$ as the base is also that $\ln(x)$ is the **area** under the curve up to $x$ (at least for $x > 1$).

So, in summary, I have shown that it is useful to express numbers as either $10^x$ or as $p^x$ (where $p$ is a positive integer other than 10). This was relatively easy to understand if $x$ is an integer because $10^x$ is a short way of writing a number that is a 'power' of 10, $i.e.$ it is $(10 \times 10 \times 10 \times ...... \times 10)_{x\,factors}$, which we usually write as a 1 with $x$ zeroes following it. For base $p$ it is a number that is a power of $p$, $i.e.$ it is $(p \times p \times p \times ...... \times p)_{x\,factors}$. I hope to have made plausible that we can extend the definition to non-integer $x$. If we name $p^x = y$ then we write $x = \log_p(y)$. When $p = 10$ we speak of *Briggsian logarithms*. When $p = e \approx 2.718$ we speak of *natural logarithms*.

A table of logarithms was handy in the pre-calculator days because, if you had to multiply two large numbers, $e.g.$ $y_1 = 3157$ by $y_2 = 24,551$, then you would look up their logarithms (base 10) and find $x_1 \approx 3.4993$ and $x_2 \approx 4.3901$. Add these up to get $7.8894$ for the sum and look what number (the 'antilog') has this as its $\log_{10}$ to find it is $10^{7.8894} = 77,507,507$. Actually this inverse lookup is done in an easier way. We know

that $7.8894 = 7 + 0.8894$ so that the answer is $10^{7.8894} = 10^{7+0.8894} = 10^7 \times 10^{0.8894}$. So we only need to look up the antilog of $0.8894$ as it is easy to multiply it by $10^7$. And in general we only need to look up the antilog of numbers between $0$ and $1$ so that the log table does not have to be huge, only of the numbers between $0$ and $1$ to the desired accuracy!